\documentclass[12pt]{article}
\usepackage{amsthm}
\usepackage{amsmath}
\usepackage{amsfonts}
\usepackage{tikz}
\usepackage{fullpage}

\newtheorem{thm}{Theorem}

\newtheorem{con}[thm]{Conjecture}

\begin{document}
\title{Conjectured lower bound for the clique number of a graph}
\author{Clive Elphick\thanks{clive.elphick@gmail.com, School of Mathematics, University of Birmingham, Birmingham, UK} \quad Pawel Wocjan\thanks{wocjan@cs.ucf.edu, Department of Computer Science, University of Central Florida, Florida, USA}}

\date{April 10, 2018}

\maketitle

\abstract{It is well known that $n/(n - \mu)$, where $\mu$ is the spectral radius of  a graph with $n$ vertices, is a lower bound for the clique number. We conjecture that $\mu$ can be replaced in this bound with $\sqrt{s^+}$, where $s^+$ is the sum of the squares of the positive eigenvalues. We prove this conjecture for various classes of graphs, including triangle-free graphs, and for almost all graphs.}

\section{Introduction}

Let $G$ be  a graph, with no isolated vertices, with $n$ vertices, edge set $E$ with $|E| = m$, average degree $d$, chromatic number $\chi(G)$ and clique number $\omega(G)$. 
We also let $A$ denote the adjacency matrix of $G$ and let $\mu = \mu_1 \ge \ldots \ge \mu_n$ denote the eigenvalues of $A$. The inertia of $A$ is the ordered triple $(\pi, \nu, \gamma)$ where $\pi$, $\nu$ and $\gamma$ are the numbers counting muliplicities of positive, negative and zero eigenvalues of $A$ respectively. Let

\[
s^+ = \sum_{i=1}^\pi \mu_i^2  \quad\mbox{and}\quad s^- = \sum_{i=n-\nu+1}^n \mu_i^2.
\]

Note that:

\[
\sum_{i=1}^n \mu_i^2 = \mathrm{tr}(A^2) = 2m = s^+ + s^-.
\]

\section{Replacing $\mu^2$ with $s^+$}

Edwards and Elphick \cite{edwards83} proved that

\[
\frac{2m}{2m - \mu^2} \le \chi(G)
\]

and Ando and Lin \cite{ando15} proved a conjecture due to Wocjan and Elphick \cite{wocjan13} that

\[
\frac{2m}{2m - s^+} = 1 + \frac{s^+}{s^-} \le \chi(G).
\]

As another example of replacing $\mu^2$ with $s^+$, Hong \cite{hong93} proved for graphs with no isolated vertices that $\mu^2 \le 2m - n + 1$, and Elphick \emph{et al} \cite{elphick16} proved that for almost all connected graphs $s^+ \le 2m - n + 1$. Similarly Favaron \emph{et al} \cite{favaron93} proved that $\omega(G) \le 2m/\mu$ and Wu and Elphick \cite{wu17} proved the doubly stronger result that $\chi(G) \le 2m/\sqrt{s^+}.$ Finally Stanley \cite{stanley87} proved that 

\[
\mu \le \frac{\sqrt{8m + 1} - 1}{2},
\]

and Wu and Elphick \cite{wu17} proved that

\[
\sqrt{s^+} \le \frac{\sqrt{8m + 1} - 1}{2}.
\]

So in all of these cases we can strengthen known bounds by replacing $\mu^2$ with $s^+$. The next section considers the same replacement for a well known lower bound for the clique number.

\section{Conjectured bound for the clique number}

The concise version of Tur\'an's theorem states that:

\begin{equation}\label{turan}
\frac{n}{n - d} \le  \omega(G).
\end{equation}

This bound was improved by Caro \cite{caro79} and Wei \cite{wei81} using degrees as follows:

\[
\sum_{i=1}^n \frac{1}{n - d_i} \le \omega(G);
\]

and by Wilf \cite{wilf86} using the spectral radius as follows:

\begin{equation}\label{wilf}
 \frac{n}{n - \mu} \le \omega(G).
\end{equation}

Bound (2) was strengthened by Nikiforov \cite{nikiforov02} who proved the following conjecture of Edwards and Elphick \cite{edwards83}.

\begin{equation}\label{vlado}
\frac{2m}{2m - \mu^2} \le \omega(G).
\end{equation}

Note that for regular graphs, all of these bounds equal $n/(n - d)$. Wocjan and Elphick \cite{wocjan13} noted that

\[
\frac{2m}{2m - s^+} \not\le \omega(G).
\]

An alternative strengthening of Wilf's bound is provided by the following conjecture, which we have tested against the thousands of named graphs with up to 40 vertices in the Wolfram Mathematica database, and found no counter-example. Aouchiche \cite{aouchiche16} has tested this conjecture using his powerful AGX software, and also found no counter-example. Conjecture 1 exceeds $n/(n - d)$ for all regular graphs with more than one positive eigenvalue.
\begin{con}
For any graph $G$

\[
\frac{n}{n - \sqrt{s^+}} \le \omega(G).
\]
This conjecture is exact, for example,  for complete regular multipartite graphs.

\end{con}

We can prove this conjecture for the following classes of graphs.

\subsection{Proof for triangle-free graphs}
\begin{proof}

Let $t$ denote the number of triangles in a graph. It is well known that:

\[
\sum_{i=1}^n \mu_i^3 = \mathrm{tr}(A^3) = 6t,  
\]

so for triangle-free graphs

\[
\sum_{i=1}^\pi \mu_i^3 = -\sum_{i=n-\nu+1}^n \mu_i^3.
\]

Therefore, using that $\mu \ge |\mu_n|$

\[
s^- \ge \frac{\sum_{i=n-\nu+1}^n \mu_i^3}{\mu_n} = \frac{\sum_{i=1}^\pi \mu_i^3}{|\mu_n|} \ge \frac{\mu^3}{|\mu_n|} \ge \mu^2.
\]

Therefore, using the lower bound on the largest eigenvalue $\mu \ge 2m/n$, the equality $\frac{1}{2}(s^+ + s^-) = m$ combined with the arithmetic-geometric-mean inequality, and the above lower bound on $s^-$, we obtain
\[
\sqrt{s^+} \le \mu \frac{n}{2m} \sqrt{s^+} \le \frac{n}{2m} \sqrt{s^-} \sqrt{s^+} \le \frac{n}{2m} \frac{2m}{2} = \frac{n}{2}.
\]

\end{proof}

\subsection{Proof for weakly perfect graphs}

\begin{proof}
Weakly perfect graphs have $\omega(G) = \chi(G)$. Therefore using the result due to Ando and Lin \cite{ando15} discussed above and that $\mu \ge 2m/n$:

\[
\frac{n}{n - \sqrt{s^+}} \le \frac{2m}{2m - s^+} \le \chi(G) = \omega(G).
\]

\end{proof}

\subsection{Proof for some strongly regular graphs} 

We do not know how to prove this conjecture for all strongly regular graphs. However we can prove the conjecture for the subset of strongly regular graphs which are Kneser graphs. The Kneser graph $KG_{p,k}$ is the graph whose vertices correspond to the $k-$element subset of a set of $p$ elements, and where two vertices are joined if and only if the corresponding sets are disjoint. The Kneser graphs with $k = 2$ are strongly regular, with only three distinct eigenvalues. For these graphs

\[
n = {p \choose 2},    \omega = \left\lfloor \frac{p}{2}\right \rfloor,  2m = {p \choose 2}{p - 2 \choose 2} \mbox{  and  } p \ge 2k = 4.
\]

The eigenvalues (see Godsil and Royle \cite{godsil01}) are:

\[
(-1)^i{p - 2 - i \choose 2 - i} \mbox{  with multiplicity  } {p \choose i} - {p \choose i - 1}, \mbox{  for  } i = 0,1,2.
\]

We are seeking to prove that

\[
\frac{n}{n - \sqrt{s^+}} \le \frac{p-1}{2} \le \left\lfloor \frac{p}{2}\right \rfloor = \omega(KG_{p,2}),
\]

which re-arranges to

\[
s^+ = 2m - s^- \le \frac{n^2(p - 3)^2}{(p - 1)^2} = \frac{p^2(p - 3)^2}{4}.
\]

Inserting the negative eigenvalues this becomes

\[
{p \choose 2}{p - 2 \choose 2} - (p - 1){p - 3 \choose 1}^2 \le \frac{p^2(p - 3)^2}{4}.
\]

Simple algebra reduces this to 

\[
2p^2 - 9p + 6 \ge 0
\]

which is true for all $p \ge 4$.

\subsection{Proof for almost all graphs}
\begin{proof}

We use the Erdos-Renyi random graph model $G_n(p)$, which consists of all graphs with $n$ vertices in which edges are chosen independently with probability $p$. Bollob\'as and Erdos \cite{bollobas76} proved that the clique number is almost always $x$ or $x+1$ where

\[
x = \frac{2\log n}{\log (1/p)} + O(\log \log n).
\]

Since almost all graphs have all degrees very close to $n/2$ we let $p = 0.5$. Therefore
\[
s^+ \le s^+ + s^- = 2m \approx n \, \frac{n}{2}
\]

So for almost all graphs

\[
\frac{n}{n - \sqrt{s^+}} \le \frac{n}{n - n/\sqrt{2}} \approx 3.4 <\frac{2\log n}{\log 2} \approx \omega(G).
\]

\end{proof}

\section{Conclusion}

Lower bounds for the clique number are often proved using the Motzkin-Straus \cite{motzkin65} inequality, which can be expressed as follows. For any adjacent vertices $i$ and $j$ such that $i < j$ we write $i \sim j$. Then for any vector $(p_1, \ldots, p_n)$ with $p_i \ge 0$ for all $i$ and $\sum_{i=1}^n p_i = 1$:

\[
\sum_{i \sim j} p_i p_j \le \frac{\omega - 1}{2\omega}.
\]

It is however not evident how to use this approach in the context of Conjecture 1, where the number of positive eigenvalues varies greatly between graphs with $n$ vertices.

\section*{Acknowledgements}

This research was supported in part by the National Science Foundation Award 1525943 .


\begin{thebibliography}{10}

\bibitem{ando15}
T. Ando and M. Lin, \emph{Proof of a conjectured lower bound on the chromatic number of a graph}, Lin. Algebra and Appl., 485, (2015), 480  - 484.

\bibitem{aouchiche16}
M, Aouchiche, private correspondence, (2016).

\bibitem{bollobas76}
B. Bollob\'as and P. Erdos, \emph{Cliques in random graphs}, Math. Proc. Camb. Phil. Soc., 80, (1976), 419 - 427.

\bibitem{caro79}
Y. Caro, \emph{New results on the independence number} Technical report, Tel Aviv University, (1979).

\bibitem{edwards83}
C. Edwards and C. Elphick, \emph{Lower bounds for the clique and the chromatic number of a graph}, Discrete Appl. Math, 5 (1983), 51 - 64.

\bibitem{elphick16}
C. Elphick, M. Farber, F. Goldberg and P. Wocjan, \emph{Conjectured bounds for the sum of squares of positive eigenvalues of a graph}, Discrete Math., 339, (2016), 2215 - 2223.

\bibitem{favaron93}
O. Favaron, M. Mah\'eo, J-F. Sacl\'e, \emph{Some eigenvalue properties in graphs (conjectures in Graffiti II)}, Discrete Math, 111, (1993), 197 - 220.

\bibitem{godsil01}
C. Godsil and G. Royle, \emph{Algebraic Graph Theory}, Springer-Verlag, New York, 2001.

\bibitem{hong93}
Y. Hong, \emph{Bounds on eigenvalues of graphs}, Discrete Math, 123, (1993), 65 - 74.

\bibitem{motzkin65}
T. Motzkin and E. Straus, \emph{Maxima for graphs and a new proof of a theorem of Tur\'an}, Canad. J. Math., (1965), 533 - 540.

\bibitem{nikiforov02}
V. Nikiforov, \emph{Some inequalities for the largest eigenvalue of a graph}. Combin. Probab. Comput. 11 (2002), 179 - 189.

\bibitem{stanley87}
R. P. Stanley, \emph{A bound on the spectral radius of graphs with $e$ edges}, Linear Algebra Appl., 87, (1987), 267 - 269.

\bibitem{wei81}
V. K. Wei, \emph{A lower bound on the stability number of a simple graph}, Technical Report 81-11217-9, Bell Laboratories, (1981).

\bibitem{wilf86}
H. Wilf, \emph{Spectral bounds for the clique and independence numbers of graphs}, J. Combin. Theory Ser. B, 40, (1986), 113 - 117.

\bibitem{wocjan13}
P. Wocjan and C. Elphick, \emph{New spectral bounds on the chromatic number encompassing all eigenvalues of the adjacency matrix}, Elec. J. Combinatorics, 20(3), (2013), P39.

\bibitem{wu17}
B. Wu and C. Elphick, \emph{Upper bounds for the achromatic and coloring numbers of a graph}, Discrete Appl. Math, 217, (2017), 375 - 380.


\end{thebibliography}
\end{document}